\documentstyle[12pt]{article}
\newcommand{\Ibb}[1]{ {\rm I\ifmmode\mkern
            -3.6mu\else\kern -.2em\fi#1}}
\newcommand{\ibb}[1]{\leavevmode\hbox{\kern.3em\vrule
     height 1.2ex depth -.3ex width .2pt\kern-.3em\rm#1}}

\newcommand{\be}{\begin{eqnarray}}
\newcommand{\ee}{\end{eqnarray}}
\newcommand{\bez}{\begin{eqnarray*}}
\newcommand{\eez}{\end{eqnarray*}}
\renewcommand{\O}{\Omega}
\newcommand{\A}{{\cal A}}
\renewcommand{\d}{\mbox{d}}
\newcommand{\Hom}{\mbox{Hom\,}}
\newcommand{\id}{\mbox{id\,}}
\newcommand{\End}{\mbox{End\,}}
\newcommand{\na}{\nabla}
\newcommand{\hna}{\hat{\nabla}}

\newcommand{\oc}{\otimes}
\newcommand{\oa}{\otimes_\A}
\begin{document}

\begin{tabbing}
\hspace*{12cm} \= GOET-TP 3/96  \\
                          \>March 1996
\end{tabbing}
\vskip1.cm

\centerline{\LARGE \bf A note on Connections and Bimodules} 

\vskip1cm

\begin{center}
      {\large \bf Aristophanes Dimakis}
       \vskip.3cm 
      Institut f\"ur Theoretische Physik  \\
      Bunsenstr. 9, D-37073 G\"ottingen, Germany
\end{center} 

\begin{abstract}
\noindent
We discuss three problems related to connections on bimodules. These
are left and right Leibniz rule for connections, left and right
linearity of their curvatures and extension of connections to 
tensor products of modules. 
\end{abstract}

\vskip1cm 
\noindent
{\bf 1.} Recently in a series of papers \cite{mour}, \cite{debo1}
and \cite{mado} the concept of a connection on a
bimodule was studied. If $\A$ is an algebra and $(\O(\A),\d)$ a 
differential calculus over $\A$, and $M$  an $\A$-bimodule, a right-module
connection on $M$ is a map 
\be \na:M\to M\oa\O^1(\A)\;,  \ee
satisfying for $f\in\A$ and $a\in M$ the {\em right Leibniz rule}
\be  \na(a f) = (\na a)f + a\oa\d f\;. \ee
In order to take account of the left action of $\A$ on $M$ the
authors of the above papers use a generalized permutation 
\be  \sigma:\O^1(\A)\oa M\to M\oa\O^1(\A)\;, \ee
and impose additionally a {\em left Leibniz rule}
\be \na(f a)=\sigma(\d f\oa a)+f \,\na a\;. \label{rlr}\ee
It turns out that this implies a strong restriction on $\na$,
a fact which reveals a close correlation between the differential calculus 
over $\A$ and the geometry on $M$.

Evidence for (\ref{rlr}) appeared also in \cite{bres}. There it was shown 
that under general conditions, if $N$ is a right $\A$-module with a 
connection $\na':N\to N\oa\O^1(\A)$, then a necessary and sufficient 
condition for $\na$ and $\na'$ to define a connection 
$\na_\oc:N\oa M\to N\oa M\oa\O^1(\A)$ in a fashion similar to the 
commutative case is that $\na$ satisfies the left Leibniz rule (\ref{rlr}) 
and
\be 
\na_\oc(b\oa a):= (\id_N\oa\sigma)(\na' b)\oa a + b\oa \na a\;. 
\ee

A further  point related to the left Leibniz rule is the tensor
character of the curvature of the connection. Extending	$\na$ to
\be \na:M\oa\O(\A)\to M\oa\O(\A) \;, \ee
by
\be  \na(a\oa\omega)=(\na a)\omega + a\oa\d\omega\;, \ee
the curvature of $\na$ is the map $\na^2:M\oa\O\to M\oa\O$.
It is easy to see that although this is a right $\O$-homomorphism
it is not even a left $\A$-homomorphism. In fact if one tries to
calculate left linearity, one finds that the action of $\na$ on
$\sigma$ must be specified. There is no evidence however 
how this has to be done.
This is the reason why in \cite{debo2} it is proposed the factoring out of 
the submodule of $M\oa\O^2(\A)$ generated by $\na^2(f a)- f\na^2 a$, 
for all $f\in\A$ and $a\in M$. On the quotient space the curvature is both
left and right $\A$-linear.

In the sequel we try to shed some new light on these constructions
using the fact that for an $\A$-bimodule $M$ the left action of $\A$
on $M$ induces a canonical imbedding of $\A$ into the algebra of right 
$\A$-endomorphisms of $M$.

\vskip.6cm 
\noindent
{\bf 2.} Let $\End^\A(M)$ denote the ring of right $\A$-endomorphism 
of $M$, i.e.\ for $\phi\in\End^\A(M)$, $a\in M$ and $f\in\A$ we have
$\phi(a\, f)=\phi(a)\, f$. A right module connection on $M$ induces
a map
\be   \hna :\End^\A(M) \to \Hom^\A(M,M\oa\O^1(\A))\;, \ee
defined by
\be   \hna(\phi) := \na\circ\phi - \phi\circ\na\;, \ee
where  $\phi$ in the last term above denotes 
also the obvious {\em extension} of $\phi:M\to M$ to 
$\phi:M\oa\O^1(\A)\to M\oa\O^1(\A)$ given by $\phi(a\oa\alpha):=
(\phi a)\oa\alpha$.
$\Hom^\A(M,M\oa\O^1(\A))$ denotes the additive group of
right $\A$-homorphisms.	With the aid of the extension mentioned,
it becomes in a natural way a $\End^\A(M)$-bimodule.

Since $M$ is a $\A$-bimodule there is a canonical algebra homomorphism
$\kappa_0:\A\to\End^\A(M)$ with
\be 
f\mapsto \hat{f}:=\kappa_0(f)\;,\qquad \hat{f}(a):=fa\;.
\ee
The action of $\hna$ on $\hat{f}$ is given by
\be  (\hna\hat{f})a =\na(fa)-f\na a\;. \label{hf}\ee
Let $\O_\na^1(\A)$ be the $\A$-bimodule, which is the additive subgroup of 
$\Hom^\A(M,M\oa\O^1(\A))$, generated by $(\kappa_0\A)(\hna\kappa_0\A)
(\kappa_0\A)$, and where for $f,g\in\A$ and $\Phi\in\O_\na^1(\A)$ we set 
\be f\,\Phi\, g:=\hat{f}\circ\Phi\circ \hat{g}\;. \label{f1}\ee
Let $\d_\na:=\hna\circ\kappa_0$, then it is clear that 
\be \d_\na:\A\to \O_\na^1(\A)\;, \ee
is a derivation since
\be   
  \hna(\hat{f}\circ\hat{g}) =(\hna\hat{f})\circ\hat{g}+
  \hat{f}\circ(\hna\hat{g})\;.
\ee
Hence $(\O_\na^1(\A),\d_\na)$ defines a first order differential
calculus over $\A$, 
\be \na(fa)=(\d_\na f)a +f\na a\;, \label{hnaf}\ee
and the left Leibniz rule uses this new differential calculus instead
of $\O(\A)$.

If $(\O^1_u(\A),\d_u)$ denotes the univeral 
first order differential calculus over $\A$, then there is a 
homomorphism $\kappa_1:\O^1_u(\A)\to\O_\na^1(\A_M)$ such
that $\kappa_1\circ\d_u = \d_\na$, i.e.\ the following 
diagram commutes
\be 
   \begin{array}{ccc}
\A & \buildrel {\rm d}_u \over \longrightarrow & \O^1_u(\A) \\
\parallel & &  
\kappa_1 \downarrow \phantom{\kappa_1} \\
\A & \buildrel {\rm d}_\na \over \longrightarrow & \O_\na^1(\A)
\end{array}\label{d1}
\ee
Define now $\sigma_u:\O^1_u(\A)\oa M \to M\oa\O^1_u(\A)$ by
\be	\sigma_u(\alpha\oa a):=\kappa_1(\alpha)(a)\;,  \ee
for $\alpha\in\O^1_u(\A)$ nd $a\in M$, then (\ref{hnaf}) takes the form
\be \na(fa)=\sigma_u(\d_u f\oa a) + f\na a\;. \ee
The existence of $\sigma_u$ was first proved in \cite{debo1}.
For $\sigma$ as in (\ref{rlr}) to exist, $\kappa_1$ must
factor uniquely through the projection $\pi_1:\O^1_u(\A)\to\O^1(\A)$,
that is $\kappa_1=\hat{\kappa}_1\circ\pi_1$ with
a homomorphism $\hat{\kappa}_1:\O^1(\A) \to \O_\na^1(\A)$ and
$\d_\na=\hat{\kappa}_1\circ\d$. Then
\be  \sigma(\alpha\oa a):=\hat{\kappa}_1(\alpha)(a)\;, \ee
for $\alpha\in\O^1(\A)$ and $a\in M$.

We conclude that a connection $\na$ on a bimodule $M$ satisfies a left 
Leibniz rule (\ref{rlr}), if $(\O^1_\na(\A),\d_\na)
\preceq(\O^1(\A),\d)$. Here $\preceq$ is the obvious partial ordering 
on the set of classes of isomorphic first order differential calculi 
over $\A$, where $(\O^1_1(\A),\d_1)\preceq(\O^1_2(\A),\d_2)$ if  
a homomorphism $\rho:\O^1_2(\A)\to\O^1_1(\A)$ exists, such that 
$\d_1=\rho\circ\d_2$.

\vskip.6cm 
\noindent
{\bf 3.} As $\na:M\to M\oa\O^1(\A)$ extends to $\na:M\oa\O(\A)\to
M\oa\O(\A)$ so also $\hna:\End^\A(M)\to\Hom^\A(M,M\oa\O^1(\A))$
extends to 
\be  \hna:\End^\O(M\oa\O(\A))\to\End^\O(M\oa\O(\A))\;, \ee
where $\End^\O(M\oa\O(\A))$ denotes the ring of right 
$\O(\A)$-homomorphisms of $M\oa\O(\A)$. Denote by
$\End_r^\O(M\oa\O(\A))$ the additive group of 
endomorphisms\footnote{That is for $\Phi\in\End_r^\O(M\oa\O(\A))$ and
for all $\xi\in M\oa\O^s(\A)$ with arbitrary $s$ we have 
$\Phi(\xi)\in M\oa\O^{r+s}(\A)$.} 
of degree $r$, then $\End^\O(M\oa\O(\A))$ is a graded ring and 
$\End_0^\O(M\oa\O(\A))\cong\End^\A(M)$ by the extension 
of $\phi\in\End^\A(M)$ to $\phi\in\End_0^\O(M\oa\O(\A))$ with
\be   \phi(a\oa\omega):=\phi(a)\oa\omega\;. \ee
For $\Phi\in\End^\O(M\oa\O(\A))$ of degree $r$ we set
\be  \hna\Phi := \na\circ\Phi - (-1)^r\Phi\circ\na\;. \ee
Now let $\hat{\O}$ denote the graded subalgebra of $\End^\O(M\oa\O(\A))$ 
generated by $\kappa_0\A$ and $\hna$, i.e.\ the elements of 
$\hat{\O}$ are finite linear combinations of expressions of the form
\bez 
(\hna^{k_1}\hat{f}_1)\circ \cdots\circ (\hna^{k_r}\hat{f}_r)\;,
\eez
with $\hna^0\hat{f}:=\hat{f}$.	 Obviously $\hna$ restricts
to $\hna:\hat{\O}\to\hat{\O}$. Note that although $\hna$ is a graded 
derivation, i.e.\ for $\Phi\in\hat{\O}^r$ and $\Psi\in\hat{\O}$ we have
\be  
\hna(\Phi\circ\Psi) = (\hna\Phi)\circ\Psi +
(-1)^r\Phi\circ(\hna\Psi)\;,
\ee
the pair $(\hat{\O},\hna)$ fails to be a differential algebra
over $\kappa_0\A$ since 
\be (\hna^2\hat{f})(a) = \na^2(fa) - f\na^2(a)\;, \ee
does not vanish in general. 

Set $J=\bigoplus_r J^r$ for the graded right $\O(\A)$-submodule of 
$M\oa\O(\A)$ generated by elements of the form 
\bez
(\hna^2\Phi)\xi\qquad \mbox{for all}\qquad \xi\in M\oa\O(\A)\;,
\quad \Phi\in\hat{\O}\;. 
\eez
It is important to note that $\na\,J\subset J$ and $\Phi\,J\subset J$
for all $\Phi\in\hat{\O}$. Hence taking the factor module 
$\O(M):=[M\oa\O(\A)]/J$ with canonical projection $p$, the elements  
$\Phi\in\hat{\O}$ have unique factorizations $p\circ\Phi=
\hat{\Phi}\circ p$ with 
\bez \hat{\Phi}:\O(M)\to \O(M)\;. \eez
The same holds for $\na$, although we do not introduce a new symbol
for $\na:\O(M) \to \O(M)$. Note that $\O(M)$
is a left $\A$- right $\O(\A)$-bimodule.

Let $\O_\na(\A)$ be the graded algebra of the factors
$\hat{\Phi}$ for all $\Phi\in\hat{\O}$. Note that since
$J^0=J^1=\{0\}$, we put
$\O^0_\na(\A)=\A$; $\O^1_\na(\A)$ is the $\A$-bimodule defined in
(\ref{f1}). We denote with $\hat{p}:\hat{\O} \to \O_\na(\A)$ the 
algebra homomorphism which sends $\Phi$ to $\hat{\Phi}$. It is easy 
to see now that $\hna$ factors through $\hat{p}$. We set 
$\hat{p}\circ \hna = \d_\na \circ \hat{p}$ and find 
that $\d_\na$ is a graded derivation of $\hat{\O}(\A_M)$ which satisfies 
additionlly $\d_\na{}^2=0$. Hence the pair $(\O_\na(\A),\d_\na)$
defines a differential calculus over $\A$ and the diagram in 
(\ref{d1}) extends to 
\be 
   \begin{array}{ccc}
\O_u(\A) & \buildrel {\rm d}_u \over \longrightarrow & \O_u(\A) \\
\kappa \downarrow \phantom{\kappa} & &  
\kappa \downarrow \phantom{\kappa} \\
\O_\na(\A) & \buildrel {\rm d}_\na \over \longrightarrow & \O_\na(\A)
\end{array}
\ee
Now let 
\bez \sigma_u: \O_u(\A)\oa \O(M) \to \O(M)\;, \eez
be the right $\O(\A)$-homomorphism defined for $\omega\in\O_u(\A),\;
\xi\in \O(M)$ by
\be   \sigma_u(\omega\oa\xi):= \kappa(\omega)(\xi)\;. \ee
Then $\sigma_u$ is graded in the first factor, right $\O(\A)$-linear in
the second factor with $\sigma_u(f\oa\xi)=f\xi$, and satisfies 
for $\omega_1,\omega_2\in\O_u(\A)$ and $\omega\in\O_u^r(\A)$
\be 
\sigma_u(\omega_1\omega_2\oa\xi) & = & \sigma_u(\omega_1\oa\sigma_u
(\omega_2\oa\xi)) \\
\na\sigma_u(\omega\oa\xi) & = & \sigma_u(\d_u\omega\oa\xi) + 
(-1)^r\sigma_u(\omega\oa\na\xi)\;.
\ee
Again if $\kappa=\hat{\kappa}\circ\pi$ holds with $\hat{\kappa}:\O(\A)\to
\O_\na(\A)$, then for $\omega\in\O(\A)$
\be   \sigma(\omega\oa\xi):= \hat{\kappa}(\omega)(\xi)\;, \ee
defines a map $\sigma: \O(\A)\oa \O(M) \to \O(M)$
with similar properties as $\sigma_u$.

Here again we find that a $\sigma$ exists if, and only if  we have
$(\O_\na(\A),\d_\na)\preceq (\O(\A),\d)$, where now $\preceq$ is a 
partial ordering on the set of all classes of isomorphic differential 
calculi over $\A$, defined in a similar way as $\preceq$ in the last 
paragraph for first order differential calculi. In fact this set 
has the structure of a lattice, induced on it by the lattice of
differential ideals of $(\O_u(\A),\d_u)$.  

\vskip.6cm 
\noindent
{\bf 4.} Let $N$ be a right $\A$-module and let
\bez \na'_M:N\to N\oa\O^1_\na(\A)\;, \eez
be a connection on $N$ with respect to the differential calculus over
$\A$ defined by $\na:M\to M\oa\O^1(\A)$. Let  $N_0:=\{b\in N|\,
b\oa a =0\;\:\mbox{for all}\;\:a\in M\}$ and $M_0:=\{a\in M|\,
b\oa a =0\;\:\mbox{for all}\;\:b\in N\}$, in the following we assume
that $\na M_0\subset M_0\oa\O^1(\A)$ and 
$\na'_M N_0\subset N_0\oa\O_\na^1(\A)$.

A connection on the tensor product $N\oa M$
\bez \na_\oc:N\oa M \to N\oa M\oa\O^1(\A)\;,\eez
is defined by
\be 
\na_\oc(b\oa a):=(\na'_M b)a + b\oa\na a\;, 
\ee
where the first term on the right side has to be interpreted in the
sense of $(b\oa\hat{\Phi})a:=b\oa(\hat{\Phi}a)$, for $\hat{\Phi}\in
\O^1_\na(\A)$.

As explained before, in case 
$(\O_\na(\A),\d_\na)\preceq(\O(\A),\d)$ we have a homomorphism 
$\hat{\kappa}:\O(\A)\to\O_\na(\A)$ with $\d_\na=\hat{\kappa}\circ\d$,
then we have also a homomorphism $\hat{\nu}:N\oa\O(\A)\to 
N\oa\O_\na(\A)$ with 
\be  \hat{\nu}(b\oa\omega) = b\oa\hat{\kappa}(\omega)\;,\ee
for $\omega\in\O(\A)$.

If $\na':N\to N\oa\O^1(\A)$ is a connection on $N$, with respect
to the original differential calculus over $\A$ then one can 
define a connection $\na_\oc:N\oa M\to N\oa M\oa\O(\A)$ 
by setting
\be
\na_\oc(b\oa a):=\hat{\nu}(\na'b)\,a + b\oa\na a\;,
\ee
if $\na'N_0\subset N_0\oa\O^1(\A)$.
In this case $\na'$ defines a unique associated connection $\na'_M$ 
making the diagram
\be 
   \begin{array}{ccc}
N\oa\O(\A) & \buildrel \na' \over \longrightarrow & N\oa\O(\A) \\
\hat{\nu} \downarrow \phantom{\hat{\nu}} & &  
\hat{\nu} \downarrow \phantom{\hat{\nu}} \\
N\oa\O_\na(\A) & \buildrel \na'_M \over \longrightarrow 
& N\oa\O_\na(\A)
\end{array}
\ee
commute.  

\vskip.6cm 
\noindent
{\bf 5.} We found that a right-module connection $\na$ on a bimodule $M$ 
over an algebra $\A$, with differential calculus $(\O(\A),\d)$ defines a new
differential calculus $(\O_\na(\A),\d_\na)$ over $\A$, in general different 
from the original calculus. With respect to this calculus $\na$ satisfies
a left Leibniz rule. If there is an algbera homomorphism
$\rho:\O(\A)\to\O_\na(\A)$, such that $\d_\na\circ\rho=\rho\circ\d$ then
we can use $\rho$ to interprete the left Leibniz rule to be with respect to 
$(\O(\A),\d)$. Otherwise we must live with two or more different 
differential calculi on $\A$ simultaneously.
In any case it is important to study, how far the induced
differential calculus $(\O_\na(\A),\d_\na)$ can be from the original
one $(\O(\A),\d)$ for given bimodule $M$.

\vskip.6cm
\noindent
{\bf Acknowledgment.} I would like to thank F.~M\"uller-Hoissen 
for stimulating discussions and  critical comments.

\end{document}